\documentclass{article}
\usepackage{graphicx}
\usepackage{pstricks}
\usepackage{float}
\usepackage{amssymb}
\usepackage{amsthm}
\usepackage{wrapfig}
\usepackage{url}
\usepackage{hyperref}

\newtheorem*{Z_stab}{Conjecture}

\title{Numerical Stability and Catalan Numbers}

\author{Arash Ghasemi\footnote{Doctorate Candidate}, Kidambi Sreenivas\footnote{Research Professor} and Lafayette K. Taylor\footnote{Professor} \\ SimCenter: National Center for Computational Engineering \\
701 E. M.L. King Blvd. Chattanooga, TN 37403, UTC}
\date{}

\begin{document}

\maketitle 

\begin{abstract}
To predict allowable time-step size for the fully discretized nonlinear differential equations, a stability theory is developed using exact determination of an infinite perturbation series. Mathematical induction is used to determine the coefficients of the series. It is discovered that the closed-form equation for the nonlinear shift of generic polynomial non-linearity can be written as a series expansion where the coefficients are the Pfaff-Fuss-Catalan numbers in Combinatorics. This reveals criteria which can be used to analytically determine the allowable time step. It is shown that stability region decreases when the nonlinearity of the differential equation increases. Therefore, the maximum allowable time step is severely limited by the nonlinearity even if an unconditionally stable scheme (in a linear sense) is used. The theory is applied to the case of general system of time-dependent nonlinear Partial Differential Equations.        
\end{abstract}

{\bf Keywords}
Numerical Stability, Nonlinear Differential Equations, Catalan Numbers, Time Marching

\section{Introduction}\label{sec:nonlin_DPI} The spatial discretization of the system of partial differential equations 
\begin{equation}
\label{first_non_PDE}
\frac{\partial v_k}{\partial t} = G\left(v_k, \frac{\partial^i v_k}{\partial x^i}, \frac{\partial^j v_k}{\partial y^j}, \ldots\right)
\end{equation}
for $\vec{v} = \left(v_1, v_2, \ldots\right) = v_k$ and some range of $i, j , \ldots$ leads to the system of semi-discrete form $d\mathbf{v}/dt = R(\mathbf{v}(t))$ or
\begin{equation}
\label{CauchyProblem}
\mathbf{v} = \mathbf{v}_0 + \int_{t_0}^{t} R\left(\mathbf{v}(\xi)\right) d\xi,
\end{equation} where $\mathbf{v} = \mathbf{v}(t)$ is the spatially distributed nodal/modal solution vector at time $t$ and $\mathbf{v}_0 = \mathbf{v}(t_0)$ designates the initial condition of the system and $R$ is the residual of the spatial discretization. The integral in (\ref{CauchyProblem}) can be arbitrarily discretized to obtain the following space-time discretization 
\begin{equation}
\label{CauchyProblem_SP}
\mathbf{u} = \mathbf{u}_0 + \Delta t \mathbb{S} \otimes R\left( \mathbf{u} \right)
\end{equation}
where $\mathbf{u} = \left[\mathbf{v}(t_1)^{T}, \mathbf{v}(t_2)^{T}, \ldots, \mathbf{v}(t_s)^{T}\right]^{T}$ is the space-time vector containing solution $\mathbf{v}(t_i)$ at temporal collocation point $t_0 < t_i \le t$, $\mathbf{u}_0 = \left[\mathbf{v}_0^{T}, \mathbf{v}_0^{T}, \ldots, \mathbf{v}_0^{T}\right]^{T}$ is the space-time initial condition, $\Delta t = (t- t_0)/s$ is the time-step, and $\mathbb{S}$ is the integration operator. For the first-order truncated Riemannian integration, $\mathbb{S}$ is a lower-diagonal unity matrix while it can be a full matrix when orthogonal polynomials such as Chebyshev polynomials are used \cite{Trefethen2001}.

Equation (\ref{CauchyProblem_SP}) represents a nonlinear system of equations which requires an iterative method to be solved in practice. A \textit{sequential iterative} solution of (\ref{CauchyProblem_SP}) constitutes the Discrete Picard Iteration (DPI) which is the exact numerical counterpart of original Picard iterations to find the fixed-point of nonlinear system (\ref{CauchyProblem_SP}). Thus, 
\begin{equation}
\label{expl_DPI}
\mathbf{u}_{n+1} = \mathbf{u}_0 + \Delta t \mathbb{S} \otimes R\left( \mathbf{u}_{n} \right)
\end{equation}
is considered here as the basic target system in which the stability of the iterative procedure is sought. However, a \textit{simultaneous} update to (\ref{CauchyProblem_SP}), i.e.   
\begin{equation}
\label{impl_DPI}
\mathbf{u}_{n+1} = \mathbf{u}_0 + \Delta t \mathbb{S} \otimes R\left( \bar{a}\,\mathbf{u}_{n+1} + \bar{b}\, \mathbf{u}_{n} \right),\;\;\;\;\; \bar{a}+\bar{b} = 1,
\end{equation}
yields an implicit form of the Discrete Picard Iteration which is also studied here. In particular, one is interested to know:

\begin{enumerate}
\item For the case where $R\left( \mathbf{u}_{n} \right)$ has polynomial nonlinearity, under what conditions are the iterative forms (\ref{expl_DPI}) and (\ref{impl_DPI}) stable?
\item If $\mathbb{S}$ is chosen such that (\ref{impl_DPI}) is linearly unconditional stable and assuming that the Jacobian of linearization is computed exactly, then does this imply nonlinear stability for arbitrary $\Delta t$?
\end{enumerate}

Answers to the above questions may improve the understanding of nonlinear instability of numerical methods which is important for researchers in the field of Computational Sciences. In general, the nonlinear systems (\ref{expl_DPI}) and (\ref{impl_DPI}) can be written as $F_h(u_h) = 0$ as an approximation to the original nonlinear system obtained without the discretization of the differential operators, i.e. $F(u) = 0$. Keller \cite{Keller75} used this notation to obtain the stability criteria based on Lipschitz continuous linearization. Later Lopez-Marcos et. al. \cite{Lopez88} worked on the same approach to interpret nonlinear stability based on local linear stability near the exact solution of the nonlinear system. As pointed out by Pirovino \cite{Pirovino91}, these linearization approaches have a disadvantage that the Lipschitz constant of the derivative $F_h(u_h) = 0$ must be known which is not possible in practice. To overcome this, Pirovino used the linearization approach in a neighborhood of $u_h$ to determine nonlinear stability. The approaches just mentioned here use norm-based inequalities to investigate the contraction of the nonlinear operator and corresponding stability. These inequality relations estimate upperbound for the solution behavior but not the \textit{exact} nonlinear mechanism which induces instability. Therefore the exact nonlinear shift (``thresholds'' according to \cite{Lopez88}) in the stability region and in the solution remains unanswered. The exact mechanism of generation of the nonlinear shift is important. Such knowledge might stimulate the design of faster algorithms with less stringent stability limits.

The structure of the paper is summarized as follows. In \autoref{sec_perturb_param} a loosely coupled form of (\ref{expl_DPI}) is considered where the perturbation parameter $\epsilon$ is introduced. This is a special case of the general theory presented at the end of the paper in \autoref{sec_gen_stab_PDEs_numeric}. Then the perturbation analysis is performed in \autoref{sec_ptb_analysis} and the exact nonlinear shift is obtained. In \autoref{sec_gen_arbit_poly_non}, the results of \autoref{sec_ptb_analysis} are generalized to arbitrary polynomial nonlinearity. For implicit discretization (\ref{impl_DPI}), the perturbation analysis is performed in \autoref{sec:nonlin_imp_DPI}. The main result of the paper is presented in \autoref{sec_gen_stab_PDEs_numeric} where the stability of the general time-dependent PDE (\ref{first_non_PDE}) is related to the concepts developed in (\autoref{sec_perturb_param}- \autoref{sec:nonlin_imp_DPI}).

\section{The perturbation parameter} \label{sec_perturb_param} 

Toward the stability analysis of (\ref{expl_DPI}) and (\ref{impl_DPI}) it is insightful to assume that the system is lossely coupled meaning that the k$^{th}$ state variable at the n$^{th}$ Picard iteration, i.e. $u_{k,n}(t)$ is almost independent of other variables $u_{l,n}(t), l \neq k$ in the solution vector $\mathbf{u}_{n}$. This assumption is exact when the residual arises from the discretization of an ordinary differential equations. However, this still remains as an approximation for the case of lower-order spatial discretization of PDEs where the system is very similar to a lossely coupled one. \footnote{In \autoref{sec_gen_stab_PDEs_numeric}, it will be shown that the results can be consistently extended to the more general cases (\ref{expl_DPI}) and (\ref{impl_DPI}) without such assumption}. In this case (\ref{expl_DPI}) can be written as $u_{k,n+1} = u_{k,0} + \Delta t \mathbb{S} \otimes R\left( u_{k,n} \right)$ or in compact form
\begin{equation}
\label{expl_DPI_2}
u_{n+1} = u_{0} + \Delta t \mathbb{S} \otimes R\left( u_{n} \right)
\end{equation}
where $u_{n+1}$ is a scalar. The second assumption in this section is that the number of temporal collocation points is limited to one. Consequently the integration operator $\mathbb{S} = \mathbb{S}_{1\times 1} = \lambda$ reduces to a scalar value which yields (\ref{expl_DPI_2}) to reduce to the following scalar equation

\begin{equation}
\label{expl_DPI_3}
u_{n+1} = u_{0} + \Delta t \lambda R\left( u_{n} \right).
\end{equation}        
Assuming that the residual $R$ is analytic over the time span, (\ref{expl_DPI_3}) can be expanded as
\begin{equation}
\label{res_anal_expanded2}
R(u_n)= R(u_0) + \left.\frac{\partial R}{\partial u}\right|_{u_0} {\Delta u}_n + \frac{1}{2}\left.\frac{\partial^2 R}{\partial u^2}\right|_{u_0} {\Delta u}_n^2 + \ldots
\end{equation}
where ${\Delta u}_n = u_n - u_0$. Substituting (\ref{res_anal_expanded2}) into (\ref{expl_DPI_3}) results in
\begin{equation}
\label{time_evol_disc_20}
{\Delta u}_{n+1} =  \Delta t \: \lambda \: \left(R_0 + \left.\frac{\partial R}{\partial u}\right|_{u_0} {\Delta u}_n + \frac{1}{2}\left.\frac{\partial^2 R}{\partial u^2}\right|_{u_0} {\Delta u}_n^2 + \ldots \right).
\end{equation}

The perturbation parameter $\epsilon$ is introduced here as a relation between first derivative (Jacobian) and higher-order derivative (Hessian)\footnote{The generalization of (\ref{perturb_param_epsil}) is presented in (\ref{gen_perturb_param_epsil}).}. 

\begin{equation}
\label{perturb_param_epsil}
\left.\frac{\partial^2 R}{\partial u^2}\right|_{u_0} = 2\: \epsilon\:  \left.\frac{\partial R}{\partial u}\right|_{u_0} 
\end{equation} 

The intuition for selecting $\epsilon$ as the perturbation parameter is described as follows. One can propose that the value of $\epsilon$ should be small for a weakly nonlinear residual where the second derivative is small compared to the first derivative. To validate this proposition, consider the scenario where the nonlinear residual converges to the linear functional $R \to c\,u$ where $c$ is a constant. Then $\left.\partial R/\partial u\right|_{u_0} \to c$ and $\left. \partial^2 R/ \partial u^2\right|_{u_0} \to 0$ which means that $\epsilon \to 0$ must hold in eq.(\ref{perturb_param_epsil}) as $R \to c\,u$. However it should be noted that the analysis presented in the following sections is valid for arbitrarily large $\epsilon$ since the perturbation series is not truncated.
Substituting (\ref{perturb_param_epsil}) into (\ref{time_evol_disc_20}) yields
\begin{equation}
\label{time_evol_disc_2}
{\Delta u}_{n+1} =  \Delta t \:\lambda \: \left(R_0 + \left.\frac{\partial R}{\partial u}\right|_{u_0} {\Delta u}_n +  \epsilon \left.\frac{\partial R}{\partial u}\right|_{u_0} {\Delta u}_n^2 + \ldots \right)
\end{equation}
Defining linear stability number as
\begin{equation}
\label{def_lin_stab_num}
r = \Delta t \: \lambda \: \left.\frac{\partial R}{\partial u}\right|_{u_0}
\end{equation}
Equation (\ref{time_evol_disc_2}) can be written as
\begin{equation}
\label{time_evol_disc_3}
{\Delta u}_{n+1} =  \Delta t \:\lambda \: R_0 + r {\Delta u}_n +  \epsilon r {\Delta u}_n^2 + \ldots
\end{equation}
To simplify notation define ${U}_0 = \Delta t \:\lambda\: R_0$ and $U = \Delta u$. Hence (\ref{time_evol_disc_3}) yields 

\begin{equation}
\label{time_evol_disc_5}
{U}_{n+1} =  {U}_0 + r (1+ \epsilon {U}_n) {U}_n + \ldots
\end{equation}

This is the final form which will be analyzed using the formal perturbation technique. Note that in this case, the nonlinear residual is

\begin{equation}
\label{non_resi_scal_for_perturbation}
R_n = c U_n + c \epsilon U_n^2, 
\end{equation}
where $c = \left.\frac{\partial R}{\partial u}\right|_{u_0}$ is the Jacobian of the linearization.
 
\section{Perturbation Analysis}\label{sec_ptb_analysis}
The solution to eq. (\ref{time_evol_disc_5}) is expanded in the term of $\epsilon$ and the i$^{th}$ perturbation amplitudes at the n$^{th}$ Picard iteration, i.e., $u_{i,n}$ such that

\begin{eqnarray}
\label{u_perturb_expansion}
U_n = \sum_{i=0}^{\infty} u_{i,n} {\epsilon}^i, 
\end{eqnarray}  
subject to initial condition 
\begin{equation}
\label{u_init_cond}
U_{n=0} = u_0 + 0 \epsilon+ 0 {\epsilon}^{2} + \ldots
\end{equation}  

Substituting (\ref{u_perturb_expansion}) into (\ref{time_evol_disc_5}) and matching the coefficients of $\epsilon^i$, a cascade of linear equations is obtained which are recursively solved to find perturbation amplitudes. It is shown in Appendix (\ref{app_Der_Perturb_Amp_For_Expl_DPI}) that the i$^{th}$ perturbation amplitude converges to
\begin{equation}
\label{solution_uin_inf} 
\frac{u_{i,\infty}}{u_0^{i+1}} = C(i)\,{\frac {{{r}}^{i}}{ \left( 1- {r} \right) ^{2i+1}}} \:\:\:\: \left(i = 0,1,2,\ldots, \:\:\: \left|r \right| \leq 1\right),
\end{equation}
where $C(i) = \{1,1,2,5,14,42,132,429\ldots\}$ is the well-known Catalan sequence \cite{OEIS} given explicitly as
\begin{equation}
\label{CatlanSeq}
C(i) = \frac{(2i)!}{i! \times (i+1)!} = \frac{\textrm{binomial}\left(2 i,i\right)}{i+1}
\end{equation} 

According to \cite{OEIS}, this sequence has many different interpretations in Combinatorics but nothing about nonlinear stability of time-stepping methods has been reported so far. Substituting the perturbation amplitudes (\ref{solution_uin_inf}) into (\ref{u_perturb_expansion}), the final nonlinear solution to DPI (\ref{time_evol_disc_5}) is obtained as follows.

\begin{eqnarray}
\nonumber \frac{U}{u_0} = \sum_{i=0}^{\infty} C(i)\,{\frac {{{\it r}}^{i}}{ \left( 1- {\it r} \right) ^{2i+1}}} {\left(\epsilon u_0\right)}^i 
\end{eqnarray}
or 
  
\begin{eqnarray}
\label{u_perturb_expansion_sol}
\frac{U}{u_0} = \underbrace{\frac{1}{1-r}}_{\textrm{Linear}} + \underbrace{\sum_{i=1}^{\infty} C(i)\,{\frac {{{\it r}}^{i}}{ \left( 1- {\it r} \right) ^{2i+1}}} {\left(\hat{\epsilon}\right)}^i}_{\textrm{Nonlinear Shift}}, \;\;\;\;\; \hat{\epsilon} = \epsilon u_0 
\end{eqnarray}
where $\hat{\epsilon}$ is introduced as the \textit{combined perturbation amplitude}. Although all perturbation amplitudes converge in the linear (original) stability region $|r| < 1$ \footnote{as shown in (\ref{solution_u0n_inf}), (\ref{solution_u1n_inf}), (\ref{solution_u2n_inf}), (\ref{solution_u3n_inf}), (\ref{solution_u4n_inf}), (\ref{solution_u5n_inf}) and (\ref{solution_uin_inf})}, their partial sum identified as {\it nonlinear shift} in (\ref{u_perturb_expansion_sol}) may or may not converge in this region. Therefore one can conclude that the linear stability region is affected as a consequence of the existence of the nonlinear shift. 

In fact (\ref{time_evol_disc_5}) is stable for some stability number $r$, if the nonlinear shift in (\ref{u_perturb_expansion_sol}) remains finite for the given perturbation amplitude $\epsilon$ and initial condition $u_0$. In order to derive an exact analytical relation for the stability region, the nonlinear shift in (\ref{u_perturb_expansion_sol}) is rearranged as follows. 
 
\begin{eqnarray}
\label{u_perturb_expansion_shift_1}
\textrm{Nonlinear Shift} = \sum_{i=1}^{\infty} C(i)\,{\frac {{{\it r}}^{i}}{ \left( 1- {\it r} \right) ^{2i+1}}} {\left(\hat{\epsilon}\right)}^i = \frac {1}{ 1- {\it r} } \sum_{i=1}^{\infty} C(i)\,{\frac {{{\it r}}^{i}}{ \left( 1- {\it r} \right) ^{2i}}} {\left(\hat{\epsilon}\right)}^i
\end{eqnarray}
Substituting the Catalan sequence from (\ref{CatlanSeq}) into (\ref{u_perturb_expansion_shift_1}) yields

\begin{eqnarray}
\label{u_perturb_expansion_shift_2}
\textrm{Nonlinear Shift} =  \frac {1}{ 1- {\it r} } \sum_{i=1}^{\infty} \frac{(2i)!}{i! \times (i+1)!} \,{\left({\frac {{\it r} \hat{\epsilon}}{ \left( 1- {\it r} \right) ^{2}}} \right)}^i
\end{eqnarray}
Therefore in order to find criteria for convergence, it is only required to find the convergence of (\ref{u_perturb_expansion_shift_2}). To achieve more compact notation define
\begin{equation}
\label{theta_non_DPI}
\theta = {\frac {{\it r} \hat{\epsilon}}{ \left( 1- {\it r} \right) ^{2}}}.
\end{equation}
where $\theta$ is named here as \textit{Nonlinear Stability Number}\footnote{According to analysis in  \autoref{sec:nonlin_imp_DPI}}.  Therefore 
\begin{eqnarray}
\label{u_perturb_expansion_shift_3}
\textrm{Nonlinear Shift} =  \frac {1}{ 1- {\it r} } \sum_{i=1}^{\infty} \frac{(2i)!}{i! \times (i+1)!} \,{\left( \theta \right)}^i
\end{eqnarray}
Thus the primary goal is to find the conditions for which the above series converges. Using the generalized hypergeometric function, it can be shown that  

\begin{equation}
\label{u_perturb_expansion_shift_4}
\sum_{i=1}^{k} \frac{(2i)!}{i! \times (i+1)!} \,{\left( \theta \right)}^i = {\frac {4\,\theta}{ \left( 1+\sqrt {1-4\,\theta} \right) ^{2}}}-{
\frac { {{}_2F_1(1,k+\frac{3}{2};\,k+3;\,4\,\theta)} \,{\theta}^{k+1} \left( 2\,k+2 \right) !
}{ \left( k+1 \right) !\, \left( k+2 \right) !}}
\end{equation}
where the standard hypergeometric function \cite{Abramowitz65, Bell2004} is expanded in terms of Gamma functions as follows  
\begin{equation}
\label{geom_series}
{}_2F_1(1,k+\frac{3}{2};\,k+3;\,4\,\theta) = \sum_{j=0}^{\infty}\frac{{\left(4\theta\right)}^j \times \; \frac{\Gamma(1+j)}{\Gamma(1)}\times \frac{\Gamma(k+3/2+j)}{\Gamma(k+3/2)} }{j! \times \; \frac{\Gamma(k+3+j)}{\Gamma(k+3)}} = \sum_{j=0}^{\infty}\frac{{\left(4\theta\right)}^j \times \; {(1)}^j\times {(k+3/2)}^j}{j! \times \; {(k+3)}^j} 
\end{equation}
Each Gamma ratio is a Pochhammer symbol. Since $k\to \infty$ then
\begin{equation}
\label{geom_series_lap1}
{}_2F_1(1,k+\frac{3}{2};\,k+3;\,4\,\theta) = \sum_{j=0}^{\infty}\frac{{\left(4\theta\right)}^j \times \; {(1)}^j\times {(k+3/2)}^j}{j! \times \; {(k+3)}^j} = \sum_{j=0}^{\infty}\frac{{\left(4\theta\right)}^j}{j!} \times \; {\left(\frac{k+3/2}{k+3}  \right)}^j = \sum_{j=0}^{\infty}\frac{{\left(4\theta\right)}^j}{j!}, 
\end{equation}
which converges to
\begin{equation}
\label{geom_series_lap2}
{}_2F_1(1,k+\frac{3}{2};\,k+3;\,4\,\theta) = \exp\left(4 \, \theta\right).
\end{equation}
Substituting (\ref{geom_series_lap2}) in (\ref{u_perturb_expansion_shift_4}) yields
\begin{equation}
\label{u_perturb_expansion_shift_4_qq}
\sum_{i=1}^{k} \frac{(2i)!}{i! \times (i+1)!} \,{\left( \theta \right)}^i = {\frac {4\,\theta}{ \left( 1+\sqrt {1-4\,\theta} \right) ^{2}}}-
\frac { \exp(4\; \theta) \,{\theta}^{k+1} \left( 2\,k+2 \right) !
}{ \left( k+1 \right) !\, \left( k+2 \right) !}
\end{equation}
For $\theta \leq 1/4$, (\ref{u_perturb_expansion_shift_4_qq}) yields real values. Hence for $k \to \infty$, (\ref{u_perturb_expansion_shift_4_qq}) reduces to
\begin{equation}
\label{u_perturb_expansion_shift_5}
\sum_{i=1}^{k} \frac{(2i)!}{i! \times (i+1)!} \,{\left( \theta \right)}^i = {\frac {4\,\theta}{ \left( 1+\sqrt {1-4\,\theta} \right) ^{2}}}, \:\:\: \theta \leq \frac{1}{4},
\end{equation}
Substituting (\ref{u_perturb_expansion_shift_5}) into (\ref{u_perturb_expansion_shift_3}) the exact nonlinear shift can be written as 
\begin{eqnarray}
\label{u_perturb_expansion_shift_6}
\textrm{Nonlinear Shift} =   {\frac {4 \theta}{\left(1- {\it r}\right) \left( 1+\sqrt {1-4\,\theta} \right) ^{2}}}
\end{eqnarray}
Also the exact converged nonlinear solution is obtained by substituting (\ref{u_perturb_expansion_shift_6}) into (\ref{u_perturb_expansion_sol}). The final result is 
\begin{eqnarray}
\label{u_perturb_expansion_sol_converged}
\frac{U}{u_0} = \frac{1}{\left(1- {\it r}\right)} \left(1 + {\frac {4 \theta}{\left( 1+\sqrt {1-4\,\theta} \right) ^{2}}} \right). 
\end{eqnarray}
For the fully linear case $\epsilon = 0$ hence $\hat{\epsilon} = \epsilon u_0 = 0$ for any initial condition and therefore $\theta = 0$ which according to (\ref{u_perturb_expansion_sol_converged}), solution for the linear case is retrieved as follows.
\begin{eqnarray}
\label{u_perturb_expansion_sol_converged_2}
\frac{U}{u_0} = \frac{1}{\left(1- {\it r}\right)}  
\end{eqnarray}
However for the fully nonlinear case  $\epsilon = 1$ hence $\hat{\epsilon} = u_0$ and therefore
\begin{eqnarray}
\label{u_perturb_expansion_sol_converged_3}
\frac{U}{u_0} = \underbrace{\left(\frac{2}{1+\sqrt {{\frac {1-2\,r+{r}^{2}-4\,r{\it u_0}}{ \left( -1
+r \right) ^{2}}}}}\right)}_{\textrm{Correction Factor}} \frac{1}{1-r}
\end{eqnarray}
Note that the correction factor converges to unity as $r \to 0$ which is consistent with the fact that for the small values of the stability number (intuitively small $\Delta t$), the problem is essentially linear. Using the definition of $\theta$ in (\ref{theta_non_DPI}) one can find the stability borders as follows. Solving (\ref{theta_non_DPI}) for $r$ yields

\begin{equation}
\label{theta_relation_roots}
r_{1,2} = 1+\frac {\hat{\epsilon}}{2\theta}\pm\sqrt {\frac{\hat{\epsilon}}{\theta}+{\left(\frac{\hat{\epsilon}}{2\theta}\right)}^{2}}
\end{equation}
\begin{figure}[H]
\begin{center}
\includegraphics[trim = 1mm 1mm 0mm 5mm, clip, width=.6\textwidth]{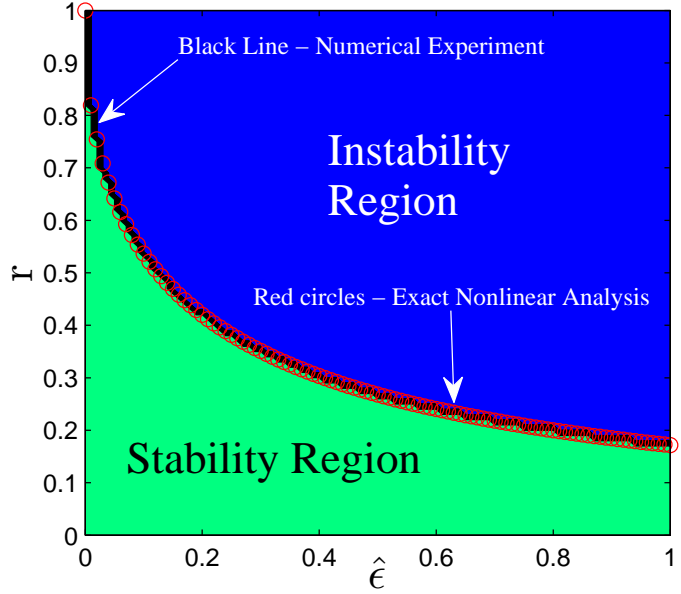}
\end{center}
\caption{The stability region of nonlinear explicit DPI (\ref{time_evol_disc_5}). This green area corresponds to (\ref{theta_relation_roots_2}). }
\label{nonlin_shift_stability_explicit_Pic}
\end{figure}

For positive perturbations $\hat{\epsilon}$ and $\theta$, the first root $r_1$ corresponding to the plus sign in (\ref{theta_relation_roots}) violates $|r| < 1$, i.e., the convergence interval of (\ref{solution_u0n_inf}, \ref{solution_u1n_inf}, \ref{solution_u2n_inf}, \ref{solution_u3n_inf}, \ref{solution_u4n_inf}, \ref{solution_u5n_inf}, \ref{solution_uin_inf}). Therefore only the second root is acceptable. Hence

\begin{equation}
\label{theta_relation_roots_2}
r \leq 1+{\frac {\hat{\epsilon}}{2\theta}}-\sqrt {\frac{\hat{\epsilon}}{\theta}+{\left(\frac{\hat{\epsilon}}{2\theta}\right)}^{2}}
\end{equation}    
The above equation determines the stability region which is plotted in green in fig.(\ref{nonlin_shift_stability_explicit_Pic}). For small perturbation amplitude $\epsilon \to 0$ and/or small initial condition $u_0 \to 0$, the combined perturbation amplitude $\hat{\epsilon} = \epsilon \, u_0  \to 0$ and therefore the linear stability condition $r \leq 1$ is retrieved by vertical axis $(\hat{\epsilon} = 0, r)$ according to fig.(\ref{nonlin_shift_stability_explicit_Pic}). The border between stability and instability regions is obtained by substituting $\theta = 1/4$ into (\ref{theta_relation_roots_2}) which yields
 
\begin{equation}
\label{nonlin_explicit_dpi_instab_border}
r = 1+2\,\hat{\epsilon}-2\,\sqrt {\hat{\epsilon}+{\hat{\epsilon}}^{2}}
\end{equation}

The result is a parabola and is plotted in fig.(\ref{nonlin_shift_stability_explicit_Pic}) using red circles. This is in exact agreement with the values (black line) obtained from the numerical solution to (\ref{time_evol_disc_5}) using a brute-force method for parameters $0 \le r,\hat{\epsilon} \le 1$.  

\section{Generalization to polynomial nonlinearity} \label{sec_gen_arbit_poly_non}
The stability analysis of nonlinear explicit DPI presented in the previous section can be consistently extended to the more general case where the residual is assumed to be a polynomial function of the dependent variable. The result is presented as follows.

\begin{Z_stab}
\label{the_conjucture}
The exact solution to the following explicit Discrete Picard Iteration
\begin{equation}
\label{time_evol_disc_f}
{U}_{n+1} =  {U}_0 + r \left(1+ \epsilon {U}_n^{Z}\right) {U}_n, \:\:\: Z = 1, 2, 3, \ldots
\end{equation}
is
\begin{eqnarray}
\label{u_perturb_expansion_sol_conj}
\frac{U}{u_0} = \left(1  + \sum_{i=1}^{\infty} C(i,Z)\,{\theta}^i\right) \frac{1}{1-r},
\end{eqnarray}    
where $\hat{\epsilon} = \epsilon u_0^{Z} $ is the combined perturbation amplitude and 
\begin{equation}
\label{gen_theta}
\theta = \frac {r \hat{\epsilon} }{{\left( 1-r \right)}^{Z+1}}
\end{equation}
is the Nonlinear Stability Number and $C(i,Z)$ is a generalized form of Catalan sequence given as
\begin{equation}
\label{gen_CatlanSeq}
C(i,Z) = \frac{\textrm{binomial}\left((Z+1) \times i,i\right)}{Z \times i+1} = \frac{((Z+1) \times i)!}{i! \times \left(Z\times i+1\right)!}.
\end{equation} 
In addition, the stability border is the solution to
\begin{equation}
\label{final_gen_Z_stab_border_canonical}
{\tilde{r}}^{Z+1} + b \tilde{r} = b, 
\end{equation} 
where $\tilde{r} = 1 - r$ and $b = \frac{{\left( Z+1 \right)}^{{Z+1}}}{Z^Z} \hat{\epsilon}$.
\end{Z_stab}

The above conjecture is validated for $Z = 1,2,3$ using symbolic processing \cite{TN_maple}. The generalized Catalan sequence given in (\ref{gen_CatlanSeq}) is known in Combinatorics as the Pfaff-Fuss-Catalan or k-Raney sequence \cite{Lang2009}. It is used in Graph Theory to enumerate (Z-ary) trees (rooted, ordered, incomplete) with $Z$ vertices including the root \cite{OEIS}.

This conjecture shed light on the mechanisms of nonlinear numerical instability. Obviously, the stability is governed by the convergence of $\sum_{i=1}^{\infty} C(i,Z)\,{\theta}^i$ in (\ref{u_perturb_expansion_sol_conj}). To understand this, it is better to find the converged value of the series for $Z = 2,3,4,5,\ldots$ by method of mathematical induction. Case $Z=1$ was studied before. For $Z=2$, one can write
\begin{equation}
\label{C_iZ_series_Z_2}
\sum_{i=1}^{\infty} C(i,2)\,{\theta}^i = \theta \times
{\mbox{$_3$F$_2$}(1,\frac{4}{3},\frac{5}{3};\,2,\frac{5}{2};\,{\frac {27}{4}}\,\theta)}
\end{equation}
where the above hypergeometric series \mbox{$_3$F$_2$} is convergent if ${\frac {27}{4}}\,\theta \leq 1$. Therefore the stability border is obtained as
\begin{equation}
\label{stab_dorder_Z_2}
{\theta}_{max}(2) = \frac{4}{27} = \frac{Z^Z}{{\left(Z+1\right)}^{Z+1}}.
\end{equation} Similarly $Z=3$ yields 
\begin{equation}
\label{C_iZ_series_Z_3}
\sum_{i=1}^{\infty} C(i,3)\,{\theta}^i = \theta \times
{\mbox{$_4$F$_3$}(1,\frac{5}{4},\frac{3}{2},\frac{7}{4};\,\frac{5}{3},2,\frac{7}{3};\,{\frac {256}{27}}\,\theta)}
\end{equation} which is convergent for 
\begin{equation}
\label{stab_dorder_Z_3}
{\theta}_{max}(3) = \frac{27}{256} = \frac{Z^Z}{{\left(Z+1\right)}^{Z+1}}.
\end{equation}
For $Z= 4$ the partial sum reduces to 
\begin{equation}
\label{C_iZ_series_Z_4}
\sum_{i=1}^{\infty} C(i,4)\,{\theta}^i = \theta \times
{\mbox{$_5$F$_4$}(1,\frac{6}{5},\frac{7}{5},\frac{8}{5},\frac{9}{5};\,\frac{3}{2},\frac{7}{4},2,\frac{9}{4};\,{\frac {3125}{256}}\,\theta)}
\end{equation}
which yields 
\begin{equation}
\label{stab_dorder_Z_4}
{\theta}_{max}(4) = \frac{256}{3125} = \frac{Z^Z}{{\left(Z+1\right)}^{Z+1}}.
\end{equation}
Similarly for $Z= 5$ one obtains 
\begin{equation}
\label{C_iZ_series_Z_5}
\sum_{i=1}^{\infty} C(i,5)\,{\theta}^i = \theta \times
{\mbox{$_6$F$_5$}(1,\frac{7}{6},\frac{4}{3},\frac{3}{2},\frac{5}{3},{\frac {11}{6}};\,\frac{7}{5},\frac{8}{5},\frac{9}{5},2,{\frac {11}{5}};\,{\frac {46656}{3125}}\,\theta)}
\end{equation}
which is convergent for 
\begin{equation}
\label{stab_dorder_Z_5}
{\theta}_{max}(5) = \frac{3125}{46656} = \frac{Z^Z}{{\left(Z+1\right)}^{Z+1}},
\end{equation} Therefore it is concluded that for arbitrary $Z$ 
\begin{equation}
\label{stab_dorder_Z_Z}
{\theta}_{max}(Z) = \frac{Z^Z}{{\left(Z+1\right)}^{Z+1}}.
\end{equation}

\begin{figure}[H]
\begin{center}
\includegraphics[trim = 3mm 2mm 10mm 6mm, clip, width=.4\textwidth]{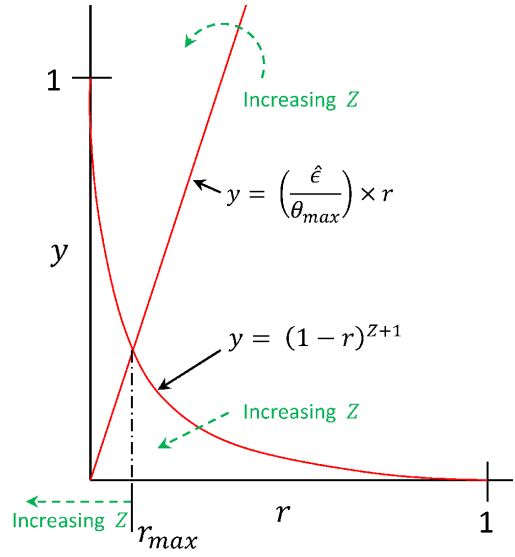}
\end{center}
\caption{The effect of increasing the degree of nonlinearity `$Z$' on the maximum allowable stability number. }
\label{r_max_shift_Z_change}
\end{figure}

To understand the effect of \textit{increasing nonlinearity}, i.e. $Z$ on the stability region, a geometrical interpretation of (\ref{gen_theta}) is possible. At stability border $\theta = {\theta}_{max}$ or
\begin{equation}
\label{gen_theta_bord}
{\theta}_{max} = \frac {r \hat{\epsilon} }{{\left( 1-r \right)}^{Z+1}}
\end{equation}
where ${\theta}_{max}$ is given in (\ref{stab_dorder_Z_Z}). Once the above equation is solved the maximum allowable stability number $r_{max}$ can be precisely determined. Unfortunately (\ref{gen_theta_bord}) is $Z+1$ degree polynomial equation and can't be solved analytically. However geometrical interpretative tools can be used. Here (\ref{gen_theta_bord}) is rearranged to define function $y$ as the below 
\begin{equation}
\label{gen_theta_bord2}
 y = {\left( 1-r \right)}^{Z+1} = \left( \frac{\hat{\epsilon} }{{\theta}_{max} } \right) r
\end{equation}
The set of curves $y = {\left( 1-r \right)}^{Z+1}$ and $y= (\hat{\epsilon}/{\theta}_{max})r$ intersect at some point $0\le r_{max}<1$ which is a solution to the original unsolvable nonlinear equation (\ref{gen_theta_bord}). This is schematically shown in fig.(\ref{r_max_shift_Z_change}) where the red curves represent the lhs and rhs of (\ref{gen_theta_bord2}).

According to fig.(\ref{r_max_shift_Z_change}), with increasing $Z$ the value of `${\theta}_{max}$' defined in (\ref{stab_dorder_Z_Z}) decreases, hence the slope of the straight line increases. On the other hand, increasing $Z$ forces the curve to ``bow'' closer to the origin. As the total result, the point of intersection of these two curves moves closer to the origin and this proves that the maximum stability number decreases.     

It should be noted that the stability border is always a \textit{canonical curve}. This can be easily shown by writing

\begin{equation}
\label{gen_theta_border}
\frac {r \hat{\epsilon} }{{\left( 1-r \right)}^{Z+1}} = {\theta}_{max} = \frac{Z^Z}{{\left(Z+1\right)}^{Z+1}} 
\end{equation}
or
\begin{equation}
\label{final_gen_Z_stab_border}
{\left(\frac{1-r}{1+Z}\right)}^{1+Z} - \left(\frac{\hat{\epsilon}}{Z^Z}\right) r = 0,
\end{equation} where the stability region is specified by 
\begin{equation}
\label{final_gen_Z_stab_border_xx}
r \leq \textrm{RootsOf} \left[ {\left(\frac{1-r}{1+Z}\right)}^{1+Z} = \left(\frac{\hat{\epsilon}}{Z^Z}\right) r\right].
\end{equation} Equation (\ref{final_gen_Z_stab_border}) can be written in the standard canonical form by changing variable $\tilde{r} = 1 - r$

\begin{equation}
\label{final_gen_Z_stab_border_canonical}
{\tilde{r}}^{Z+1} + b \tilde{r} = b, 
\end{equation} where $b$ is a constant given as
\begin{equation}
\label{final_gen_Z_stab_border_canonical_b}
b = \frac{{\left( Z+1 \right)}^{{Z+1}}}{Z^Z} \hat{\epsilon}.
\end{equation}

\section{Nonlinear Stability Analysis of Implicit DPI} \label{sec:nonlin_imp_DPI} 
The perturbation analysis of \autoref{sec_ptb_analysis} can be applied to the case where DPI is performed implicitly. In this case the residual vector in (\ref{impl_DPI}) can be written as the weighted average between two iterative steps
\begin{equation}
\label{impl_DPI2}
\mathbf{u}_{n+1} = \mathbf{u}_0 + \Delta t \mathbb{S} \otimes \left(a\,R( \mathbf{u}_{n+1} )+ b \, R( \mathbf{u}_{n} ) \right).
\end{equation}
where weights satisfy $a+b = 1$ and $0\le a,b \le 1$. For $a = 0$ and $b = 1$ the explicit DPI (\ref{expl_DPI}) is retrieved. For $a = 1$ and $b = 0$, (\ref{impl_DPI2}) yields
\begin{equation}
\label{impl_DPI3}
\mathbf{u}_{n+1} = \mathbf{u}_0 + \Delta t \mathbb{S} \otimes R( \mathbf{u}_{n+1} ).
\end{equation}
Substituting the linearization of residual, i.e. $R( \mathbf{u}_{n+1} ) \simeq R( \mathbf{u}_{n} ) + \partial R( \mathbf{u}_{n} )/ \partial \mathbf{u}_{n} \left(\mathbf{u}_{n+1} - \mathbf{u}_{n} \right)$ in (\ref{impl_DPI3}) yields
\begin{equation}
\label{impl_DPI4}
\left(\mathbb{I} - \Delta t \mathbb{S} \otimes \left.\frac{\partial R}{ \partial \mathbf{u}_{n}}\right|_{\mathbf{u}_{n}} \right) \mathbf{u}_{n+1} = \mathbf{u}_0 + \Delta t \mathbb{S} \otimes \left(R( \mathbf{u}_{n} ) - \left.\frac{\partial R}{ \partial \mathbf{u}_{n}}\right|_{\mathbf{u}_{n}} \mathbf{u}_{n} \right)
\end{equation}
Following th assumptions made in \autoref{sec_perturb_param}, (\ref{impl_DPI4}) reduces to

\begin{eqnarray}
\label{time_evol_disc_gen_compact_imp_linear_final_stab_analysis}
\bigg(1 - \Delta t\:\lambda \: \frac{\partial R_{n}}{\partial U_n} \bigg) U_{n+1} = u_0 + \Delta t \: \lambda \: \left(R_{n} - \frac{\partial R_{n}}{\partial U_n} U_{n}\right)
\end{eqnarray}
where the Jacobian is the derivative of the residual defined in eq.(\ref{non_resi_scal_for_perturbation})

\begin{equation}
\label{jacob_non_resi_scal_for_perturbation}
\frac{\partial R_n}{\partial U_n} = c + 2 \;c \;\epsilon\; U_n = c \left(1+2\;\epsilon\; U_n \right), 
\end{equation}
Substituting (\ref{jacob_non_resi_scal_for_perturbation}) and (\ref{non_resi_scal_for_perturbation}) into (\ref{time_evol_disc_gen_compact_imp_linear_final_stab_analysis}) results in

\begin{eqnarray}
\label{time_evol_disc_gen_compact_imp_linear_final_stab_analysis_2}
\bigg(1 - \Delta t\:\lambda \: c \left(1+2\;\epsilon\; U_n \right) \bigg) U_{n+1} = u_0 + \Delta t \: \lambda \: \left( c U_n + c \epsilon U_n^2 - c \left(1+2\;\epsilon\; U_n \right) U_{n}\right)
\end{eqnarray}
Using the definition of the stability number $r = c \lambda \Delta t$ (\ref{time_evol_disc_gen_compact_imp_linear_final_stab_analysis_2}) simplifies to 

\begin{eqnarray}
\label{time_evol_disc_gen_compact_imp_linear_final_stab_analysis_3}
\nonumber
\bigg(1 - r \left(1+2\;\epsilon\; U_n \right) \bigg) U_{n+1} &=& u_0 - r \epsilon U_n^2 \\
U_{n = 0} &=& u_0,
\end{eqnarray}

This is the sequence that is analyzed here. The perturbation series (\ref{u_perturb_expansion}) is then substituted into (\ref{time_evol_disc_gen_compact_imp_linear_final_stab_analysis_3}) which generates expressions for perturbation amplitudes. It is shown in Appendix (\ref{app_pert_amp_imp_dpi}) that the i$^{th}$ perturbation amplitude of implicit DPI is written as
\begin{equation}
\label{uin_imp_fin} 
\frac{u_{i,n}}{u_0^{i+1}} = C(i)\,{\frac {{{r}}^{i}}{ \left( 1- {r} \right) ^{2i+1}}} \:\:\:\: i = 0,1,2,\ldots,
\end{equation}
Where $C(i)$ are the Catalan numbers. Comparing the above with (\ref{solution_uin_inf}) it is clear that both results are equal except $r$ is not constrained in (\ref{uin_imp_fin}) since its convergence is independent of $n$. Substituting (\ref{uin_imp_fin}) into the perturbation series (\ref{u_perturb_expansion}) yields

\begin{eqnarray}
\label{u_perturb_expansion_sol_implicit}
\frac{U}{u_0} = \frac{1}{1-r} + \frac {1}{ 1- {\it r} } \sum_{i=1}^{\infty} \frac{(2i)!}{i! \times (i+1)!} \,{\theta}^i, \;\;\;\;\; \hat{\epsilon} = \epsilon u_0, \;\;\;\; \theta = {\frac {{\it r} \hat{\epsilon}}{ \left( 1- {\it r} \right) ^{2}}}
\end{eqnarray}

This is consistent with explicit DPI for $|r| \leq 1$ (see (\ref{u_perturb_expansion_shift_2})). In fact the stability and convergence of implicit DPI is only governed by the convergence of the nonlinear shift in (\ref{u_perturb_expansion_sol_implicit}) not the original stability number $r$. Therefore the nonlinear stability number $\theta$ acts like a stability number governing the nonlinear nature of the residual and this is the reasoning behind its name. 

\begin{figure}[H]
\centering
\includegraphics[trim = 15mm 0mm 14mm 6mm, clip, width=0.49\textwidth]{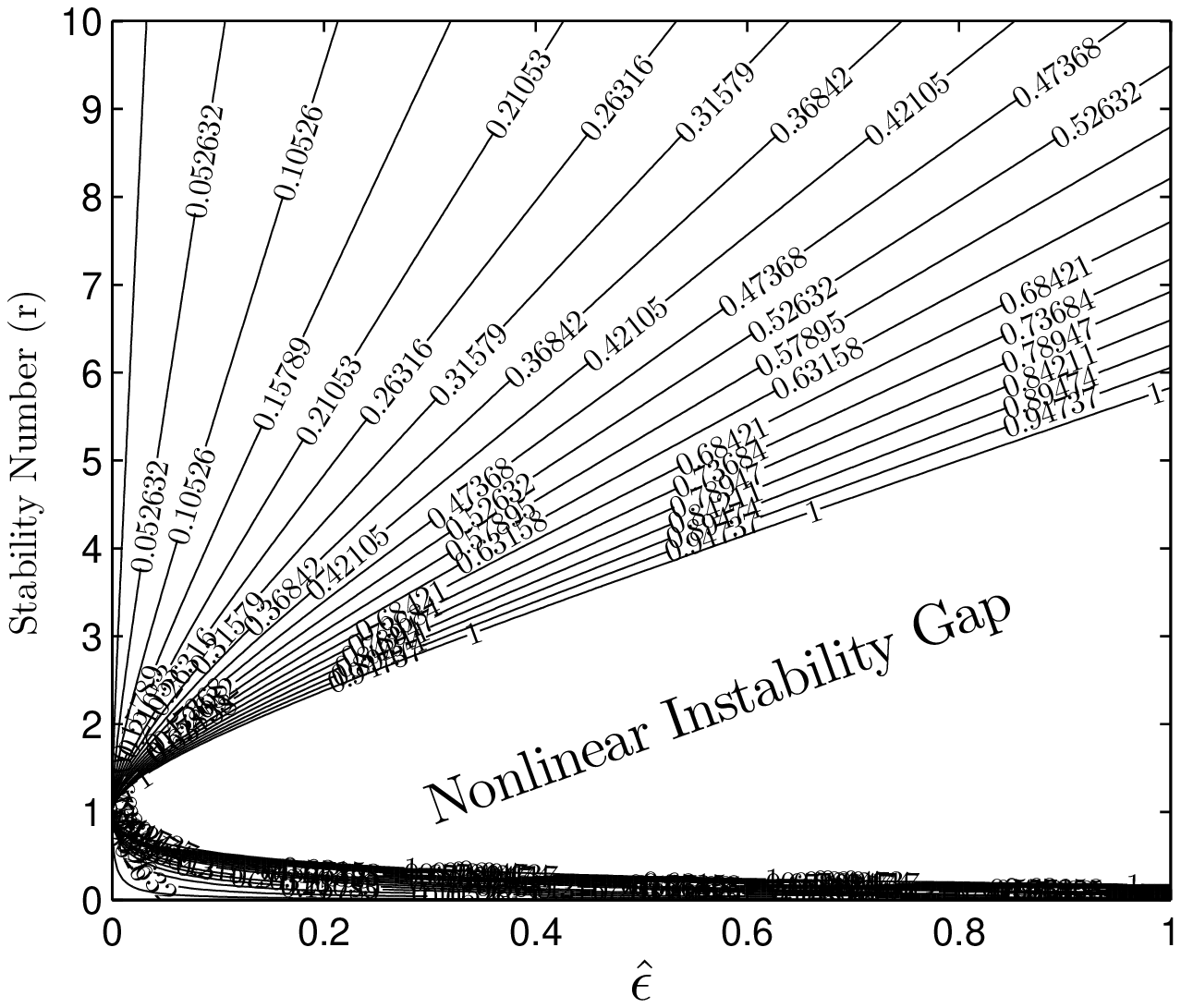}
\includegraphics[trim = 15mm 0mm 14mm 6mm, clip, width=0.49\textwidth]{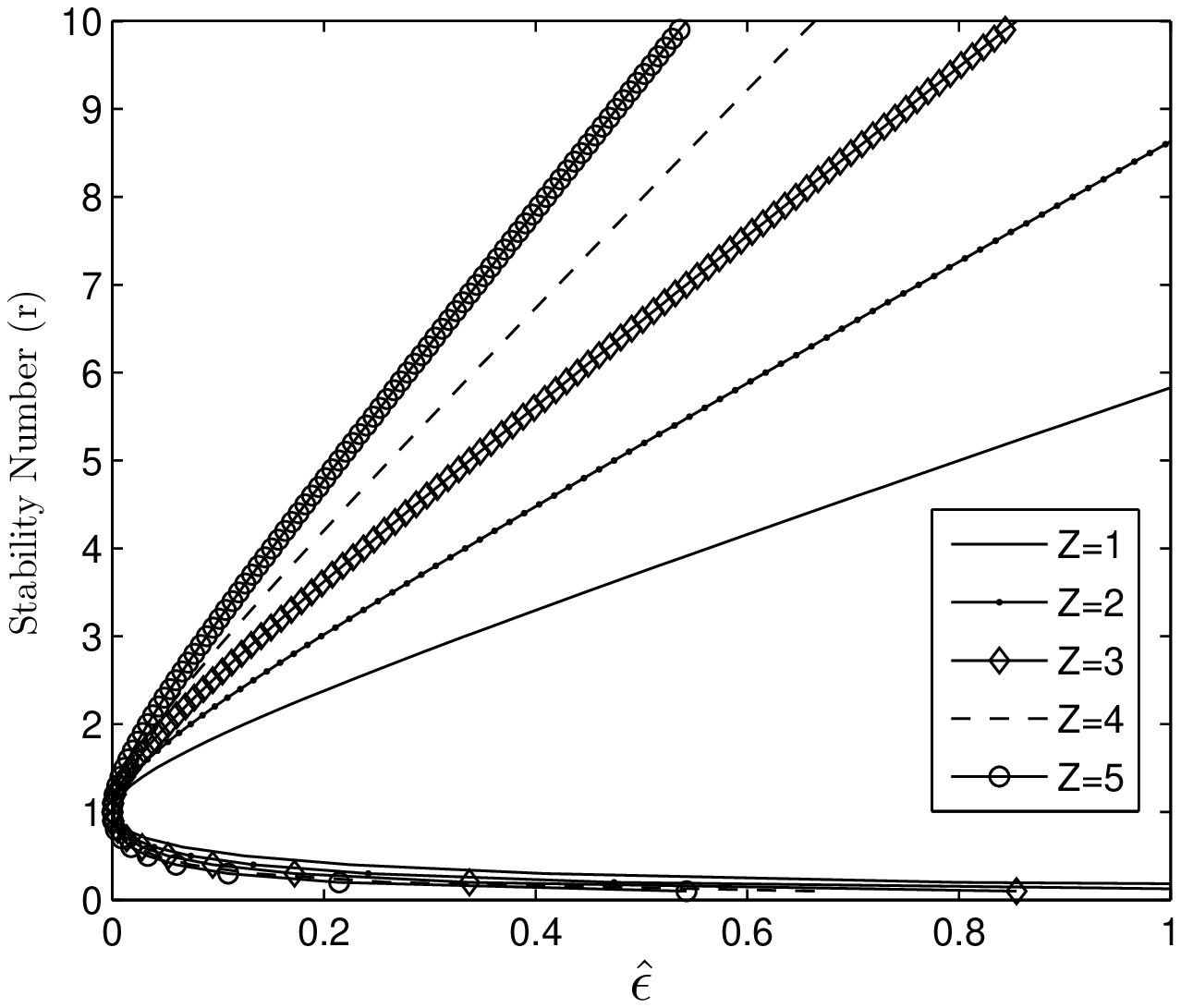}
\caption{Left) A contour plot of ${\frac {4\;{\it r} \hat{\epsilon}}{ \left( 1- {\it r} \right) ^{2}}}$ versus $\hat{\epsilon}$ and $r$. Right) Variation of maximum stability number for increasing nonlinearity. The nonlinear instability gap increases as nonlinearity $Z$ increases.}
\label{non_mod_stab_number}
\end{figure}

According to (\ref{u_perturb_expansion_shift_5}), $\theta$ must be less than or equal to 1/4 so that the nonlinear shift converges. This implies that

\begin{equation}
\label{imp_DPI_modified_non_stab_cond}
\left|{\frac {{\it r} \hat{\epsilon}}{ \left( 1- {\it r} \right) ^{2}}} \right| \leq \frac{1}{4}
\end{equation}    
The stability regions of (\ref{imp_DPI_modified_non_stab_cond}) is shown in figure (\ref{non_mod_stab_number}-Left). The area under $r=1$ is exactly equal to the nonlinear stability theory of the explicit DPI described by (\ref{theta_relation_roots_2}) and presented in fig.(\ref{nonlin_shift_stability_explicit_Pic}).  

For the case $r>1$ the nonlinear implicit DPI still remains stable. However there is a parabolic nonlinear instability gap which must be avoided in practice. This important results reveals a weakness of the implicit DPI. While it is \textit{linearly} unconditionally stable for $a=1$ and $b=0$\footnote{See fig.(\ref{non_mod_stab_number}-Left) for linear case ($\hat{\epsilon} = 0, r=\textrm{arbit.}$) is always in the stability region.}, it has a instability gap due to nonlinearity of the residuals. This analysis can be extended to general nonlinear polynomial residual where $R_n =  c\left(1+ \epsilon {U}_n^{Z}\right) {U}_n, \:\:\: Z = 1, 2, 3, \ldots$. According to eq.(\ref{stab_dorder_Z_Z}), in this case the bound for the modified nonlinear stability number $\theta$ is 

\begin{equation}
\label{stab_dorder_Z_Z_2}
{\theta}_{max}= {\frac {{\it r_{max}} \hat{\epsilon}}{ \left( 1- {\it r_{max}} \right) ^{2}}} = \frac{Z^Z}{{\left(Z+1\right)}^{Z+1}},
\end{equation}
where $Z=1$ is the second-order nonlinearity (see (\ref{imp_DPI_modified_non_stab_cond}). A plot of the stability number versus $\hat{\epsilon}$ is shown in fig.(\ref{non_mod_stab_number}-Right). As nonlinearity increases, i.e. $Z$ increases the instability gap widens rapidly.

It can be conclude that in practice the governing equations should be slightly nonlinear or with small initial conditions. In this case either $Z$, $\epsilon$ or $u_0$ are small thus $\hat{\epsilon}$ is small enough to neglect the instability gap according to fig.(\ref{non_mod_stab_number}-Right).

\section{Generalization to the stability of system of PDEs} \label{sec_gen_stab_PDEs_numeric}
The assumptions made in \autoref{sec_perturb_param} transformed the general nonlinear systems (\ref{expl_DPI}) and (\ref{impl_DPI}) into scalar equations (\ref{time_evol_disc_5}) and (\ref{time_evol_disc_gen_compact_imp_linear_final_stab_analysis}) which in this case an exact analysis was possible. However this analysis still can be utilized when (\ref{expl_DPI}) and (\ref{impl_DPI}) are considered to be system of arbitrary size. For a moment, lets assume that the PDE (\ref{first_non_PDE}) is not discretized. In this case, consider the corresponding integral form of (\ref{first_non_PDE}) i.e. $v_k = v_{0k} + \int_{t_0}^{t} G\left(v_k(\xi),\frac{\partial^i v_k(\xi)}{\partial x^i},\frac{\partial^j v_k(\xi)}{\partial y^j}, \ldots\right) d\xi$ and introduce the space-time analytical operator $\tilde{G} = \int_{t_0}^{t} G$. Then one can write
\begin{equation}
\label{PDE_anl1}
v_k = v_{0k} + \tilde{G}
\end{equation} 
The analytical residual $\tilde{G}$ is now expanded $\tilde{G} = \left. \tilde{G} \right|_{v_{0k}} + \mathbf{J} \Delta v_k + \frac{1}{2} {\Delta v_k}^T \, \mathbf{H} \, \Delta v_k + H.O.T$ where $\mathbf{J} = \left. \frac{\partial \tilde{G}}{\partial v_{k}} \right|_{v_{0k}}$ is the Jacobian and $\mathbf{H} = J\left(\nabla \tilde{G}\right)\left(v_{0k}\right)$ is the Hessian matrix evaluated at $v_{0k}$ and $\Delta v_k = v_k - v_{0k}$. This is analogous to the procedure in \autoref{sec_perturb_param} for the scalar case. Doing so (\ref{PDE_anl1}) yields
\begin{equation}
\label{PDE_anl2}
v_k = v_{0k} + \left. \tilde{G} \right|_{v_{0k}} + \mathbf{J} \Delta v_k + \frac{1}{2} {\Delta v_k}^T \, \mathbf{H} \, \Delta v_k + H.O.T,
\end{equation} 
or
\begin{equation}
\label{PDE_anl3}
{\| \Delta v_k \|}_p \leq  {\| \left. \tilde{G} \right|_{v_{0k}} \|}_p + { \| \mathbf{J} \Delta v_k \|}_p + \frac{1}{2} {\| {\Delta v_k}^T  \, \mathbf{H}  \, \Delta v_k  +  H.O.T. \|}_p
\end{equation}
On the other hand consider the following linear Sturm-Liouville problem 
\begin{equation}
\label{gen_Sturm_Liouville}
\mathbf{J} \, \Delta v_{k} = \left. \frac{\partial \tilde{G}}{\partial v_{k}} \right|_{v_{0k}} \Delta v_{k} = \lambda_{k} \Delta v_{k},    
\end{equation}
which can be solved analytically for fair broad range of PDEs with prescribed boundary conditions since it is a linear equation. The supremum of the eigenvalue spectrum of (\ref{gen_Sturm_Liouville}) is denoted by $r = \max \{ |\lambda_{k}| \}$. Therefore using (\ref{gen_Sturm_Liouville}), (\ref{PDE_anl3}) can be bounded by
 \begin{equation}
\label{PDE_anl4}
{\| \Delta v_k \|}_p \leq  {\| \left. \tilde{G} \right|_{v_{0k}} \|}_p + { \| r \Delta v_k \|}_p + \frac{1}{2} {\| {\Delta v_k}^T \, \mathbf{H}  \, \Delta v_k +  H.O.T. \|}_p
\end{equation}
According to the discussion in \autoref{sec_perturb_param} the second derivative can be related to the first derivative using a perturbation parameter. As a generalization to (\ref{perturb_param_epsil}), one can write
\begin{equation}
\label{gen_perturb_param_epsil}
\frac{1}{2}\,{\| \Delta v_k^T \, \mathbf{H} \, \Delta v_k + H.O.T \|}_p = \epsilon \, {\| \mathbf{J} \, \Delta v_k\|}_p\, {\| \Delta v_k\|}_p 
\end{equation}
for some arbitrary $\epsilon$ in the \textit{entire} space-time. The parameter $\epsilon$ is small when $v_k$ is close to a linear functional according to \autoref{sec_perturb_param}. However, as mentioned before, there is no restriction on the size of $\epsilon$ since perturbation series is not truncated. Thus (\ref{gen_perturb_param_epsil}) should always hold. Substituting (\ref{gen_perturb_param_epsil}) into (\ref{PDE_anl4}) yields
\begin{equation}
\label{PDE_anl5}
{\| \Delta v_k \|}_p \leq  {\| \left. \tilde{G} \right|_{v_{0k}} \|}_p + r { \| \Delta v_k \|}_p + \epsilon \, {\| \mathbf{J} \, \Delta v_k\|}_p\, {\| \Delta v_k\|}_p
\end{equation} 
Substituting (\ref{gen_Sturm_Liouville}) in (\ref{PDE_anl5}) yields
\begin{equation}
\label{PDE_anl6}
{\| \Delta v_k \|}_p \leq  {\| \left. \tilde{G} \right|_{v_{0k}} \|}_p + r { \| \Delta v_k \|}_p + \epsilon \, {\| r \Delta v_k \|}_p\, {\| \Delta v_k\|}_p
\end{equation}
Introducing $V = {\| \Delta v_{k} \|}_p$ and $V_0 = {\| \left. \tilde{G} \right|_{v_{0k}} \|}_p$, (\ref{PDE_anl6}) can be written as
\begin{equation}
\label{PDE_anl7}
V \leq V_0 + r V + \epsilon r V^2
\end{equation}
or
\begin{equation}
\label{PDE_anl8}
V \leq V_0 + \left(1 + \epsilon V\right)\, r V 
\end{equation}
which is analogous to (\ref{time_evol_disc_5}). The iterative class $V_{n+1} \leq V_0 + \left(1 + \epsilon V_{n}\right)\, r V_{n}$ is the explicit Analytical Picard Iteration (API) for the norm of the solution satisfying the general PDE (\ref{first_non_PDE}). Since $V_{n} \ge 0$ and $r \ge 0$ the perturbation method of \autoref{sec_ptb_analysis} can be consistently used here. Perturbing $V_{n}$ in the terms of $\epsilon$ similar to (\ref{u_perturb_expansion}) and the solving the corresponding recursive sequences one will obtain

\begin{eqnarray}
\label{gen_u_perturb_expansion_sol}
\frac{V}{V_0} \leq \frac{1}{1-r} + \sum_{i=1}^{\infty} C(i)\,{\frac {{{\it r}}^{i}}{ \left( 1- {\it r} \right) ^{2i+1}}} {\hat{\epsilon}}^i, \;\;\;\;\; \hat{\epsilon} = \epsilon V_0 
\end{eqnarray}
which remains bounded if $\left| \theta = r \hat{\epsilon}/ {\left( 1- r \right)}^{2} \right|  \leq 1/4$ according to discussion in \autoref{sec_ptb_analysis}. It can be shown that the same condition applies when the API is performed implicitly.

\begin{wrapfigure}{r}{.5\textwidth}
  \centering
\resizebox{0.4\textwidth}{!}{\input{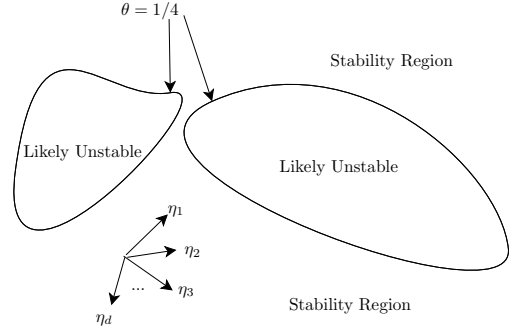}}
  \caption{The stability region of general nonlinear system of PDEs (\ref{first_non_PDE}) in the Fourier space granted by (\ref{gen_stab_PDE_wavenumber}).}
\label{fig_iso_contours}
\end{wrapfigure}
There is an interesting discussion regarding the linear Sturm-Liouville problem (\ref{gen_Sturm_Liouville}). Since $\mathbf{J}$ is always a \textit{linear operator} it can be represented via
\begin{equation}
\label{jacob_linear_v0k}
\mathbf{J} = \sum\limits_{l=1}^{d}\sum\limits_{m} a_{lm}\left(v_{0k}\right) \frac{\partial^m}{\partial x_l^m}
\end{equation}      
in the d-dimensional space $\mathbb{R}^d$. On the other hand, the Fourier transform 

\begin{equation}
\label{FTransform_in_Rd}
\hat{\Delta v_k} = \int_{\mathbb{R}^d} \Delta v_k \, e^{-\mathbf{i} x_{j} \, \eta_{j}} dx_1\,dx_2\,\ldots\,dx_d, \;\;\; j = 1\ldots d,
\end{equation}
maps $\Delta v_k$ defined in the physical domain $x_j$ to $\hat{\Delta v_k}$ in the frequency domain $\eta_j$. Taking Fourier transform of (\ref{gen_Sturm_Liouville}) and using (\ref{FTransform_in_Rd}) yields

\begin{equation}
\label{gen_Sturm_Liouville2}
\int_{\mathbb{R}^d} \mathbf{J} \, \Delta v_{k} \, e^{-\mathbf{i} x_{j} \, \eta_{j}} dx_1\,dx_2\,\ldots\,dx_d = \lambda_{k} \int_{\mathbb{R}^d} \Delta v_k \, e^{-\mathbf{i} x_{j} \, \eta_{j}} dx_1\,dx_2\,\ldots\,dx_d,
\end{equation}
or
\begin{equation}
\label{gen_Sturm_Liouville3}
\left( \sum\limits_{l=1}^{d}\sum\limits_{m} a_{lm}\left(v_{0k}\right) {\left(\mathbf{i}\,{\eta}_{l}\right)}^m \right)\,\hat{\Delta v_{k}}  = \lambda_{k}\, \hat{\Delta v_{k}}.    
\end{equation}
Therefore the k-th eigenvalue of the Sturm-Liouville problem is obtained as
\begin{equation}
\label{gen_Sturm_Liouville4}
\lambda_k = \sum\limits_{l=1}^{d}\sum\limits_{m} a_{lm}\left(v_{0k}\right) {\left(\mathbf{i}\,{\eta}_{l}\right)}^m,    
\end{equation}
and hence the stability number is $r = \max \left|\lambda_k\right|$. Therefore $\theta \leq 1/4$ implies that a solution to the system of nonlinear PDEs (\ref{first_non_PDE}) remains stable and finite in space-time $\mathbb{R}^d$ if
\begin{equation}
\label{gen_stab_PDE_wavenumber}
|\theta| = \left|\frac{\max \left|\sum\limits_{l=1}^{d}\sum\limits_{m} a_{lm}\left(v_{0k}\right) {\left(\mathbf{i}\,{\eta}_{l}\right)}^m\right| \hat{\epsilon}}{{\left(1-\max \left|\sum\limits_{l=1}^{d}\sum\limits_{m} a_{lm}\left(v_{0k}\right) {\left(\mathbf{i}\,{\eta}_{l}\right)}^m\right|\right)}^2}\right| \leq \frac{1}{4}
\end{equation}
The iso-level contours of (\ref{gen_stab_PDE_wavenumber}) for $\theta = 1/4$ (which are not necessarily closed curves) define the stability borders as depicted in fig.(\ref{fig_iso_contours}) where stability is guaranteed by (\ref{gen_stab_PDE_wavenumber}) outside of these regions. Inside these regions, however, the solution may or may not be stable because (\ref{gen_stab_PDE_wavenumber}) yields a least upperbound.

Also any numerical solution to (\ref{first_non_PDE}) is nonlinearly stable if (\ref{gen_stab_PDE_wavenumber}) is valid when the frequency $\eta_{l}$ is replaced with the frequency modified by the numerical method. Such a modified frequency can be easily obtained using Discrete Fourier Transform.

{\bf Example:} For nonlinear Poisson equation on defined on $x \in I = \left[-1,1\right]$
\begin{equation}
\label{non_poisson_eq_def}
\frac{\partial v}{\partial t} = \frac{\partial^2 R}{\partial x^2}, \;\;\; R = v + v^2,
\end{equation}
with IBVs 
\begin{equation}
\label{one_d_poiss_IBV}
v\left(x = -1, t\right) = v\left(x = 1, t\right) = 0, \;\;\; v\left(x,t=0\right) = \left(1 - x\right) \left(1 + x\right),
\end{equation}
the explicit DPI (\ref{expl_DPI}) leads to
\begin{equation}
u_{n+1} = u_{0} + \Delta t \, R\left({u}_{n}\right). 
\end{equation}
Therefore (\ref{PDE_anl1}) yields
\begin{equation}
 \tilde{G} = \Delta t \, R\left(u_{n}\right) = \Delta t \, \frac{\partial^2}{\partial x^2}\left(u_n + u_n^2\right)
\end{equation}
The Jacobian is
\begin{equation}
\label{jac_non_poisson}
J = \left. \frac{\partial \tilde{G}}{\partial u} \right|_{u_{0}} = \Delta t\, \left(\frac{\partial^2}{\partial x^2} \square + 2 \frac{\partial^2}{\partial x^2} u_0 \square \right),
\end{equation}
and hence the general Sturm-Liouville problem (\ref{gen_Sturm_Liouville}) reduces to the following

\begin{equation}
\label{strum_non_poisson}
\Delta t\, \frac{\partial^2}{\partial x^2} \left(\left(1+ 2 u_0\right) \Delta u_k\right) = \lambda_k \, \Delta u_k
\end{equation} 
with the stability criteria given by 

\begin{equation}
\label{stab_cond_non_poisson}
|\theta| = \left| \frac{r \hat{\epsilon}}{{\left(1-r\right)}^2} \right| \leq \frac{1}{4},
\end{equation}

where $r = \max|\lambda_k|$. If (\ref{strum_non_poisson}) is solved \textit{analytically} for \textit{infinite} eigenvalues $\lambda_{k=1\ldots \infty}$ then (\ref{stab_cond_non_poisson}) leads to \underline{semi-discrete stability regions}. In the semi-discrete approach, an infinite dimensional banded matrix is indeed considered for the Jacobian operator (\ref{jac_non_poisson}) and the equations are only discretized in time. This is while, in a fully discrete numerical solution, a finite-dimensional matrix (not necessarily banded) is employed. In this case, \underline{numerical stability regions} can be investigated by finding whether (\ref{stab_cond_non_poisson}) is satisfied for finite-dimensional eigen-spectrum $\lambda_{k=1\ldots M}$. These eigenvalues uniquely correspond to the numerical method used for discretization and also the type of boundary conditions used. This incorporates all details of a numerical solution in the current stability theory in a unified and consistent way. Therefore for different discretization method and/or BCs types, the discretized form of Jacobian matrix given in (\ref{jac_non_poisson}) changes and thus the eigen spectrum (\ref{strum_non_poisson}) changes and as a result, the stability regions obtained from (\ref{stab_cond_non_poisson}) changes accordingly.     

Focusing on the numerical stability, consider a symmetric second-order discretization of Laplacian $\frac{\partial^2}{ \partial x^2} \approx \frac{1}{\Delta x^2} \, \textrm{diag}\left(1,-2,1\right)_{M\times M}$ where $M+2$ collocation points (including the boundaries) are used on interval $I$. In this case (\ref{strum_non_poisson}) can be written as

\begin{equation}
\label{strum_non_poisson2}
\left[\beta \,.\, \textrm{diag}\left(1,-2,1\right) \, .\, \textrm{diag}\left(1+ 2 u_0(x_i)\right)\right] \Delta u_k = \lambda_k \, \Delta u_k
\end{equation}
where $\beta = \frac{\Delta t}{\Delta x^2}$ is the CFL number. According to \cite{wen_ch_yu}, the eigenvalues of the tridiagonal matrix $\textrm{diag}\left(1,-2,1\right)$ can be obtained as 

\begin{equation}
\label{eig_tridiag_2nd_order}
\bar{\lambda}_k = 2\left(\cos\left(\frac{k\, \pi}{M+1}\right) - 1\right) \;\;\;, k = 1\ldots M.
\end{equation}
Substituting (\ref{eig_tridiag_2nd_order}) in (\ref{strum_non_poisson2}) yields
\begin{equation}
\label{eig_spectrum_1}
\lambda_k = 2 \beta \left(\cos\left(\frac{k\, \pi}{M+1}\right) - 1\right) \left(1+2 (1-x_k) (1+x_k)\right) 
\end{equation} 
where $x_k = -1 + 2\, k / (M+1)$. For convenience, define a new variable
\begin{equation}
\label{def_gamma}
\gamma = \frac{k}{M+1}
\end{equation}  
Substituting $x_k = -1 + 2 \gamma$ and (\ref{def_gamma}) in (\ref{eig_spectrum_1}) yields
\begin{equation}
\label{eig_spectrum_2}
\lambda_k = 2 \beta \left(\cos(\gamma \pi) - 1\right) \left(1 + 8 \gamma \, (1-\gamma) \right)
\end{equation} 
Hence $r$ can be obtained by finding the maximum value of $\lambda_k$ over $I$. The extremum happens at the root of the derivative of (\ref{eig_spectrum_2}) which is a nonlinear equation. Therefore an exact solution is not possible and hence it is estimated as follows.
  
\begin{equation}
\label{value_of_r_non_poisson}
r = \max|\lambda_k| \approx 8.562 \beta
\end{equation}

At this moment, the value of $\hat{\epsilon}$ is required according to (\ref{stab_cond_non_poisson}) to complete the analysis. Since $\hat{\epsilon} = \epsilon V_0$, it is easier to compute $V_0$ and $\epsilon$ separately. The value of $V_0$ is obtained as follows.
\begin{equation}
\label{val_V0}
V_0 = {\| \left. \tilde{G} \right|_{u_0} \|}_p = {\| \Delta t \frac{\partial^2}{\partial x^2} \left(u_0 + u_0^2\right) \|}_p
\end{equation}     
For the second-order central numerical discretization used here, the Laplacian operator in (\ref{val_V0}) should be replaced with the corresponding discretized form as below.
\begin{equation}
\label{val_V0_2}
V_0 = {\| \frac{\Delta t}{\Delta x^2} \textrm{diag}(1,-2,1) \, \textrm{diag}(u_0 + u_0^2) \|}_p = \frac{\Delta t}{\Delta x^2}\, \max\left(\textrm{eig}(\textrm{diag}(1,-2,1)\, \textrm{diag}(u_0 + u_0^2))\right). 
\end{equation}
Substituting corresponding eigenvalues, (\ref{val_V0_2}) can be written as follows.
\begin{equation}
\label{val_V0_3}
V_0 = 2\beta \left| \max \left\{ \left(\cos(\gamma \pi) -1\right)\, (u_0(\gamma) + u_0^2(\gamma)) \right\}\right| 
\end{equation}
Since $x_k = -1 + 2 \gamma$ then $u_0(x_k) = u_0(\gamma) = 4 \gamma (1-\gamma)$ hence (\ref{val_V0_3}) leads to

\begin{equation}
\label{val_V0_4}
V_0 = 2\beta \left| \max \left\{ \left(\cos(\gamma \pi) -1\right)\, (4\gamma + 12 \gamma^2 - 32 \gamma^3 + 16 \gamma^4) \right\}\right| 
\end{equation}
which can be approximated as
\begin{equation}
\label{val_V0_5}
V_0 \approx 5.054 \, \beta
\end{equation}

The value of the Hessian in (\ref{gen_perturb_param_epsil}) can be obtained by taking the deivative of (\ref{jac_non_poisson}) which yields

\begin{equation}
\label{Hess_non_poisson}
H = \left. \frac{\partial J}{\partial u} \right|_{u_{0}} = 2\, \Delta t\, \frac{\partial^2}{\partial x^2} \square
\end{equation} 

Also note that $H.O.T = 0$ in (\ref{gen_perturb_param_epsil}) since higher derivatives of Hessian are identically zero. Substituting (\ref{Hess_non_poisson}) and (\ref{jac_non_poisson}) in (\ref{gen_perturb_param_epsil}) yields

\begin{equation}
\label{gen_perturb_param_epsil2}
{\| \frac{\partial^2}{\partial x^2} \square \|}_p = \epsilon \, {\| \frac{\partial^2}{\partial x^2} (\square + 2\, u_0 \square)\|}_p = \epsilon \frac{r}{\Delta t},
\end{equation}
or equivalently

\begin{equation}
\label{gen_perturb_param_epsil3}
\epsilon = \frac{\max\left| \textrm{eig}\left(\frac{\partial^2}{\partial x^2} \square\right) \right| }{ \max\left| \textrm{eig}\left(\frac{\partial^2}{\partial x^2} \left(\square + 2\,u_0 \square\right) \right) \right|} = \frac{\Delta t \max\left| \textrm{eig}\left(\frac{\partial^2}{\partial x^2} \square\right) \right| }{ r} = \frac{\Delta t \max\left| \textrm{eig}\left(\frac{\textrm{diag}(1,-2,1)}{\Delta x^2} \right) \right| }{ r}.  
\end{equation}
Hence

\begin{equation}
\label{gen_perturb_param_epsil4}
\epsilon = \frac{\beta \max\left| \textrm{eig}\left(\textrm{diag}(1,-2,1) \right) \right| }{ r}
\end{equation}
Substituting (\ref{eig_tridiag_2nd_order}) in (\ref{gen_perturb_param_epsil4}) yields
\begin{equation}
\label{gen_perturb_param_epsil5}
\epsilon = \frac{2\,\beta  \left(1-\cos(\frac{M \pi}{M+1}) \right)}{ r} \approx \frac{4 \beta}{r}
\end{equation}
Therefore from (\ref{gen_perturb_param_epsil5}) and (\ref{val_V0_5}) it is concluded that
 
\begin{equation}
\label{epsilon_hat_val}
\hat{\epsilon} = \epsilon \, V_0 = \frac{20.216 \beta^2}{r}
\end{equation}
Substituting (\ref{epsilon_hat_val}) and (\ref{value_of_r_non_poisson}) in (\ref{stab_cond_non_poisson}) yields

\begin{equation}
\label{stab_cond_non_poisson_values}
\left| \frac{20.216 \, \beta^2}{ {\left(1 - 8.562 \beta\right)}^2 } \right| \leq \frac{1}{4}
\end{equation}  
Solving (\ref{stab_cond_non_poisson_values}) it can be easily verified that the stability region is $0 \leq \beta \leq 0.0570$. This is a great reduction in the allowable CFL number $\beta$ compared to the linear Poisson equation where $0 \leq \beta \leq 0.5$. This spectacular result can not be justified using linear stability theories. 

To validate the analytical stability region $0 \leq \beta \leq 0.0570$, a computer program \cite{TN_maple} is written which solves nonlinear Poisson equation (\ref{non_poisson_eq_def}) with the given initial and boundary conditions using second-order spatial discretization. The value of CFL number is experimentally modified to find the stability region. It is found that $0 \leq \beta \leq 0.0885$ which is consistent with the analytical result since the current theory gives a least upperbound.

\section{Conclusions}
The analysis presented in this paper determines the stability region of nonlinear system of PDE (\ref{first_non_PDE}) when the corresponding space-time integral (\ref{CauchyProblem}) is discretized in explicit form (\ref{expl_DPI}) and implicit form (\ref{impl_DPI}). Important conclusions are summarized as follows.
\begin{enumerate}
\item The analysis presented in this paper determines the shift that occurs in the linear stability criteria due to the existence of nonlinear terms in residual. This shift was shown to be exact when (\ref{expl_DPI}) and (\ref{impl_DPI4}) are scalar and can be regarded as a least upperbound when (\ref{expl_DPI}) and (\ref{impl_DPI4}) are general system of equations. This answers the first question in the introduction. 
\item For both explicit and implicit discretization, there is a canonical instability gap in the $r-\hat{\epsilon}$ plane for polynomial nonlinearity (see fig.(\ref{non_mod_stab_number})-left). Outside of this region, the solution remains stable while inside of this gap, the scalar version of (\ref{expl_DPI}) and (\ref{impl_DPI4}) are guaranteed to be unstable. However the general form (\ref{expl_DPI}) and (\ref{impl_DPI4}) may or may not be unstable in this region according to fig.(\ref{fig_iso_contours}) and discussions in \autoref{sec_gen_stab_PDEs_numeric}. This result implies that even if linearization is done perfectly, and the Jacobian of linearization is computed analytically, and a linearly unconditional stable is applied for the discretization of (\ref{first_non_PDE}),  then still the resulting numerical method is nonlinearly unstable inside the instability gap. This address question (2) in the introduction.
\item The area of the instability gap increases when the degree of the nonlinearity of the residual increases (see fig.(\ref{non_mod_stab_number})-right). In this case, the space-time discretization of (\ref{first_non_PDE}) is strongly limited by nonlinear instability. However, from a practical point of view, application of a different discretization of the original Cauchy problem such as multi-step Runge-Kutta methods may or may not reduce the nonlinear instability gap. This prompts further investigation of the nonlinear instability of RK methods which may or may not be canonical.                
\end{enumerate} 
AG acknowledge help and support form SimCenter University of Tennessee at Chattanooga.

\appendix

\section{Derivation of Perturbation Amplitudes For Explicit DPI} \label{app_Der_Perturb_Amp_For_Expl_DPI}

The details of derivation of perturbation amplitudes for explicit DPI is presented as follows. For $i = 0$, the corresponding equation would be 
\begin{eqnarray}
\label{u0n}
\nonumber
u_{0,n+1} &=& r \, u_{0,n} + u_0   \\
u_{0,n=0} &=& u_0.
\end{eqnarray}
which is the only perturbation amplitude when the residual is linear, i.e. $R = c u$. Matching the coefficient of $\epsilon^1$ yields

\begin{eqnarray}
\label{u1n}
\nonumber u_{1,n+1} &=& r u_{1,n} + r u_{0,n}^{2} \\ 
u_{1,n=0}  &=& 0
\end{eqnarray}
Similarly for the coefficient of $\epsilon^2$ one obtains 
 
\begin{eqnarray}
\label{u2n}
\nonumber
u_{2,n+1} &=& r u_{2,n} + 2 \, r \, u_{0,n} \, u_{1,n} \\
u_{2, n=0} &=& 0.
\end{eqnarray}
The coefficients of $\epsilon^3$ and $\epsilon^4$ generates the following sequences. 

\begin{eqnarray}
\label{u3and4n}
\nonumber  
u_{3,n+1} &=& r\, u_{3,n} + 2\,
r\, u_{0,n} \, u_{2,n} + r \, {u_{1,n}}^{2} \\
\nonumber u_{4, n+1} &=& r\, u_{4,n} + 2\,
r\, \left(u_{2,n} \, u_{1,n} + u_{0,n} \, u_{3,n} \right)  \\
u_{3,n=0} &=& u_{4,n = 0} = 0, 
\end{eqnarray}

It should be noted that the sequences generated in this way always consist of a linear core in the form of $r\, u_{i,n}$ plus a nonlinear source term which \textit{only} depends on the previous Picard iterations. This is the desired property of the perturbation method which makes it possible to analytically obtain the i$^{th}$ nonlinear amplitude using recursive solution of linear sequences. Similar expressions can be derived for higher order terms; however, the resulting expressions are very long to be included here. A symbolic was written to derive and solve the equation for the i$^{th}$ perturbation amplitude\cite{TN_maple}.

At this point, the perturbation amplitudes need to be solved recursively. First (\ref{u0n}) is solved yielding 

\begin{equation}
\label{solution_u0n}
\frac{u_{0,n}}{u_0} = \frac {{r}^{n+1}-1 }{r-1}
\end{equation}
which is the partial sum of the first $n$ terms of geometric series obtained by recursively expanding (\ref{u0n}). Equation (\ref{solution_u0n}) converges at arbitrarily large iterations if and only if $|r| <1$. In this case the converged solution is
\begin{equation}
\label{solution_u0n_inf}
\frac{ u_{0,\infty}}{u_0} = -{\frac {1}{r-1}}
\end{equation} 
Since for the linear residual, $u_{0,n}$ is the only available perturbation amplitude it can be concluded that the sufficient linear stability requirement is $|r| <1$. The second perturbation amplitude is obtained by substituting (\ref{solution_u0n}) into (\ref{u1n}) and finding the partial sum. The final result is

\begin{equation}
\label{solution_u1n}
\frac{u_{1,n}}{u_0^2} = {\frac { \left( -2\,{r}^{2+n}+2\,{r}^{n+1} \right) n}{ \left( -1+r
 \right) ^{3}}}+{\frac {-r+{r}^{2n+2} +{r}^{n+
1}-{r}^{2+n}}{ \left( -1+r \right) ^{3}}}, 
\end{equation}          which converges to  
\begin{equation}
\label{solution_u1n_inf}
\frac{u_{1,\infty}}{u_0^2} = -{\frac {r}{ \left( r-1 \right) ^{3}}} \,\,\,\, \textrm{iff} \,\, |r| < 1
\end{equation} Substituting (\ref{solution_u1n}) and (\ref{solution_u0n}) into (\ref{u2n}) the second perturbation amplitude can be found. The final result can be written as follows.
\begin{eqnarray}
\label{solution_u2n}
\nonumber & &\frac{ u_{2,n} }{u_0^3} = -2\,{\frac { \left( -{r}^{4+n}+{r}^{3+n}+{r}^{2+n}-{r}^{n+1} \right) {
n}^{2}}{ \left( -1+r \right) ^{4} \left( {r}^{2}-1 \right) }}-2\,{
\frac { \left( 2\,{r}^{4+2n} -2\,{r}^{2+2n} -{r}^{2+n}+{r}^{4+n} \right) n}{ \left( -1
+r \right) ^{4} \left( {r}^{2}-1 \right) }} \\ & & -2\,{\frac {-{r}^{2+n}-{r}^
{4+2n} +{r}^{4+n}+{r}^{3}-{r}^{3+3n} +{r}^{2}-{r}^{3+2n} +{r}^{3+n}}
{ \left( -1+r \right) ^{4} \left( {r}^{2}-1 \right) }}  
\end{eqnarray} 
which converges to        

\begin{equation}
\label{solution_u2n_inf} 
\frac{ u_{2,\infty} }{u_0^3} = -2\,{\frac {{{r}}^{2}}{ \left( {r}-1 \right) ^{5}}} \,\,\,\, \textrm{iff} \,\, |r| < 1
\end{equation}
The partial sum of the third and higher amplitudes are exceedingly lengthy. The converged solutions are provided here. The third amplitude yields
\begin{equation}
\label{solution_u3n_inf} 
\frac{u_{3,\infty}}{u_0^4} =
-5\,{\frac {{{r}}^{3}}{ \left( {r}-1 \right) ^{7}}} \,\,\,\, \textrm{iff} \,\, |r| < 1
\end{equation}
The full partial sum of the fourth amplitude converges to
\begin{equation}
\label{solution_u4n_inf} 
\frac{u_{4,\infty}}{u_0^5} = -14\,{\frac {{{r}}^{4}}{ \left( {r}-1 \right) ^{9}}} \,\,\,\, \textrm{iff} \,\, |r| < 1
\end{equation}
and the fifth amplitude converges to 
\begin{equation}
\label{solution_u5n_inf} 
\frac{u_{5,\infty}}{u_0^6} = -42\,{\frac {{{ r}}^{5}}{ \left( { r}-1 \right) ^{11}}} \,\,\,\, \textrm{iff} \,\, |r| < 1
\end{equation}
A symbolic processor was used to derive full partial sums and finding limits where it is determined (using mathematical induction) that the i$^{th}$ perturbation amplitude converges to
\begin{equation}
\label{solution_uin_inf} 
\frac{u_{i,\infty}}{u_0^{i+1}} = C(i)\,{\frac {{{r}}^{i}}{ \left( 1- {r} \right) ^{2i+1}}} \:\:\:\: \left(i = 0,1,2,\ldots, \:\:\: \left|r \right| \leq 1\right),
\end{equation}
where $C(i) = \{1,1,2,5,14,42,132,429\ldots\}$ is the well-known Catalan sequence \cite{OEIS} given explicitly as
\begin{equation}
\label{CatlanSeq}
C(i) = \frac{(2i)!}{i! \times (i+1)!} = \frac{\textrm{binomial}\left(2 i,i\right)}{i+1}.
\end{equation} 

\section{Perturbation Amplitudes of Implicit DPI} \label{app_pert_amp_imp_dpi}

The zeroth perturbation amplitude yields
\begin{eqnarray}
\label{u0n_imp}
(1-r) u_{0,n+1} = u_0
\end{eqnarray}
which is easily solved for $u_{0,n+1}$ as
\begin{eqnarray}
\label{u0n_imp_fin}
 \frac{u_{0,n+1}}{u_0} = -\frac{1}{r-1}.
\end{eqnarray}
Evidently $u_{0,n+1}$ for implicit DPI converges in the first iteration thus it is independent of $n$. This shows that \textit{unlike} the zeroth perturbation amplitude of \textit{explicit DPI} given in (\ref{solution_u0n}), the zeroth perturbation amplitude of implicit DPI is always stable independent of the iteration number $n$. The first amplitude is obtained as
\begin{equation}
\label{u1n_imp}
u_{1,n+1} = \frac {r\,u_{0,n}  \left( -2\, u_{0, n+1} + u_{0,n}  \right) }{r-1} 
\end{equation}
Substituting (\ref{u0n_imp_fin}) into (\ref{u1n_imp}) yields
\begin{equation}
\label{u1n_imp_fin}
\frac{u_{1,n+1}}{u_0^2} = -{\frac {r}{ \left( r-1 \right) ^{3}}}  
\end{equation}
Again, by comparing (\ref{u1n_imp_fin}) with (\ref{solution_u1n}) one realizes that the first perturbation amplitude of the implicit DPI scheme is independent of the Picard iterations. In fact the first perturbation amplitude of implicit DPI for arbitrary $r$ and $n$ is \textit{exactly equal} to the first perturbation amplitude of explicit DPI when it converges (compare to (\ref{solution_u1n_inf})). 
Similarly the second perturbation amplitude is read as        
\begin{equation}
\label{u2n_imp}
u_{2,n+1} = 2\, \frac {r \left( -u_{0,n}\, u_{1,n+1} - u_{1, n}\, u_{0, n+1} +u_{0, n} u_{1, n} \right) }{r-1} 
\end{equation}
Substituting (\ref{u0n_imp_fin}, \ref{u1n_imp_fin}) into (\ref{u2n_imp}) yields
\begin{equation}
\label{u2n_imp_fin}
\frac{u_{2,n+1}}{u_0^3} = -2\,{\frac {{r}^{2}}{ \left( r-1 \right) ^{5}}}
\end{equation}
The same conclusion again holds here whereas (\ref{u2n_imp_fin}) and (\ref{solution_u2n_inf}) are equal. Similarly for the third, fourth and the fifth perturbation amplitudes are obtained as follows.
\begin{eqnarray}
\label{u3n_imp}
\nonumber
\left(\frac{r-1}{2r}\right)\,u_{3,n+1} &=& u_{0, n}\, u_{2,n}
-u_{0,n}\, u_{2, n+1}-u_{1,n} \, u_{1, n+1}
  \\
&-& u_{2,n} \, u_{0,n+1} + \frac{1}{2}\,u_{1,n}^{2}, 
\end{eqnarray}
\begin{eqnarray}
\label{u4n_imp}
\nonumber
\left(\frac{r-1}{2r}\right)\,u_{4,n+1} &=& - u_{3,n} \, u_{0, n+1} - u_{0,n} \, u_{3,n+1} - u_{1,n}\, u_{2,n+1}\\ 
& - & u_{2,n} \, u_{1, n+1} + u_{0,n}\, u_{3,n} + u_{1,n} \, u_{2,n}.
\end{eqnarray}
\begin{eqnarray}
\label{u5n_imp}
\nonumber
\left(\frac{r-1}{2r}\right)\, u_{5,n+1} &=& -u_{3,n} \, u_{1, n+1}- u_{0,n} \, u_{4,n+1} - u_{1,n} \, u_{3,n+1} \\ 
\nonumber &-& u_{2,n}\, u_{2, n+1} - u_{4,n} \, u_{0, n+1} + u_{0,n} \, u_{4,n}\\ &+& u_{1,n}\, u_{3,n} + \frac{1}{2}\, u_{2,n}^{2}. 
\end{eqnarray}
By substituting previous perturbation amplitudes into (\ref{u3n_imp}), (\ref{u4n_imp}) and (\ref{u5n_imp}) one obtains
\begin{equation}
\label{u3n_imp_fin}
\frac{u_{3,n+1}}{u_0^4} = -5\,{\frac {{r}^{3}}{ \left( r-1 \right) ^{7}}}  
\end{equation}
\begin{equation}
\label{u4n_imp_fin}
\frac{u_{4,n+1}}{u_0^5} = -14\,{\frac {{r}^{4}}{ \left( r-1 \right) ^{9}}}  
\end{equation}
\begin{equation}
\label{u5n_imp_fin}
\frac{u_{5,n+1}}{u_0^6} = -42\,{\frac {{r}^{5}}{ \left( r-1 \right) ^{11}}}  
\end{equation}
As mentioned before, the above relations are independent of iteration number and they are in fact the exact converged value of corresponding explicit DPI relations given in (\ref{solution_u3n_inf}), (\ref{solution_u4n_inf}) and (\ref{solution_u5n_inf}). In general the i$^{th}$ perturbation amplitude of implicit DPI is written as
\begin{equation}
\label{uin_imp_fin} 
\frac{u_{i,n}}{u_0^{i+1}} = C(i)\,{\frac {{{r}}^{i}}{ \left( 1- {r} \right) ^{2i+1}}} \:\:\:\: i = 0,1,2,\ldots.
\end{equation}
   

\end{document}